# Stochastic Multistage Co-optimization of Generation and Transmission Expansion Planning

Fernanda Thomé, Ricardo C. Perez, Lucas Okamura, Alessandro Soares and Silvio Binato
{fernanda, ricardo, okamura, Alessandro, silvio}@psr-inc.com.
PSR
Rio de Janeiro, Brazil

*Abstract* — This paper describes a generation and transmission (G&T) expansion planning tool based on Benders decomposition and multistage stochastic optimization: (i) A MIP-based "investment module" determines a trial expansion plan; (ii) an SDDP-based "operation module" calculates the expected operation costs for the trial plan and Benders cuts to the investment module based on the expected marginal costs of the capacity constraints at the optimal solution.

The integrated G&T planning approach is illustrated with a realistic planning study for the Bolivian power system.

## I. NOMENCLATURE

*Investment module*

| | |
|---|---|
| $x_{t,i}, x_{t,j}, x_{t,r}$ and $x_{t,k}$ | represent the construction of candidate projects (hydro $i$, thermal $j$, renewable $r$ and circuit/transformer $k$) in stage $t$. |
| $\mathbb{I}_x, \mathbb{J}_x, \mathbb{R}_x$ and $\mathbb{K}_x$. | sets of candidate projects |
| $\mathbb{I}_e, \mathbb{J}_e, \mathbb{R}_e$ and $\mathbb{K}_e$ | sets of existing generation/transmission devices |
| $m$ | Index of the investment→operation Benders decomposition iterations |
| $R^{ex}$ | set of exclusivity constraints |
| $P_c^{ex}$ | set of projects $p$ that belong to exclusivity constraint $c$ |
| $R^{as}$ | set of association constraints |
| $P_c^{as}$ | set of projects that belong to the association constraint $c$, where the decisions on projects $p1$ and $p2$ are associated |
| $R^{pre}$ | set of precedence constraints |
| $P_k^{pre}$ | set of projects that belong to precedence constraint $c$, where project $p1$ precedes project $p2$ |

*Indices*

| | |
|---|---|
| $t = 1, \ldots, T$ | time stages (typically weeks or months) |
| $\tau = 1, \ldots, T$ | intra-stage time blocks (e.g. peak/medium/low demand or 168 hours in a week) |
| $s = 1, \ldots, S$ | scenarios for each stage $t$ produced by the stochastic models (typically inflows and renewable generation; also loads, equipment availability and fuel costs) |
| $l = 1, \ldots, L$ | set of scenarios for stage $t+1$ conditioned to scenario $s$ in stage $t$ |
| $i = 1, \ldots, I$ | storage devices (typically hydro plants; also fuel storage, batteries, emission limits and some types of contracts) |
| $m \in M_i$ | set of hydro plants immediately upstream of plant $i$ |
| $j = 1, \ldots, J$ | dispatchable devices (typically, thermal plants; also, some controllable renewables and price-responsive demand) |
| $r = 1, \ldots, R$ | non-dispatchable devices (typically, wind, solar and biomass) |
| $n = 1, \ldots, N$ | transmission network buses |
| $k = 1, \ldots, K$ | transmission network components (circuits, transformers and FACTS devices such as phase shifters and smart wires) |
| $p = 1, \ldots, P$ | number of hyperplanes (Benders cuts) in the future cost function |

*Decision variables for the operation problem in stage $t$, scenario $s$*

| | |
|---|---|
| $v_{t+1,i}$ | stored volume of hydro $i$ by the end of stage $t$ |
| $u_{t,i}$ | turbined volume of hydro $i$ stage $t$ |
| $v_{t,i}$ | spilled volume of hydro $i$ in stage $t$ |
| $e_{t,\tau,i}$ | energy produced by hydro $i$ in block $\tau$, stage $t$ |
| $g_{t,\tau,j}$ | energy produced by thermal plant $j$ in block $\tau$, stage $t$ |
| $\alpha_{t+1}^l$ | present value of expected future cost from $t+1$ to $T$ conditioned to scenario $l$ in $t+1$ |

*Known values for the operation problem in stage $t$, scenario $s$*

| | |
|---|---|
| $\hat{a}_{t,i}^s$ | lateral inflow to hydro $i$ in stage $t$, scenario $s$ ($\hat{a}_t^s$ set of inflows for all hydro plants) |
| $\hat{v}_{t,i}^s$ | stored volume of hydro $i$ in the beginning of stage $t$, scenario $s$ ($\hat{v}_t^s$ set of stored volumes for all hydro plants) |
| $\overline{v}_i$ | maximum storage of hydro $i$ |
| $\overline{u}_i$ | maximum turbined outflow of hydro $i$ |
| $\phi_i$ | production coefficient ($kWh/m^3$) of hydro $i$ |
| $\overline{g}_j$ | maximum generation of thermal plant $j$ |
| $c_j$ | variable operating cost of thermal plant $j$ |
| $\hat{r}_{t,\tau,r}^s$ | energy produced by renewable plant $r$ in stage $t$, block $\tau$, scenario $s$ |
| $\hat{d}_{t,\tau}$ | demand of block $\tau$, stage $t$ |

*Multipliers*

| | |
|---|---|
| $\pi_{t,i}^h$ | multiplier of the storage balance equation of hydro $i$ (see problem formulation) |
| $\pi_{t,i}^a$ | multiplier of the conditioned inflow equation of hydro $i$ (see problem formulation) |

*pth Benders cut coefficients*

| | |
|---|---|
| $\hat{\varphi}_{t+1,i}^{hp}$ | coefficient of hydro plant $i$'s storage, $v_{t+1,i}$ |
| $\hat{\varphi}_{t+1,i}^{ap}$ | coefficient of hydro plant $i$'s inflow, $a_{t+1,i}^l$ |
| $\hat{\varphi}_{t+1}^{0p}$ | constant term |

*Stochastic streamflow model coefficients*

| | |
|---|---|
| $\hat{\mu}_{t,i}$ | mean of the lateral inflow to hydro $i$ in stage $t$. |
| $\hat{\sigma}_{t,i}$ | standard deviation of the lateral inflow to hydro $i$ in stage $t$. |
| $\hat{\rho}_{t,i}$ | serial correlation of the lateral inflow to hydro $i$ in stage $t$. |
| $\hat{\xi}_{t,i}^{l}$ | correlation matrix for the sampled residuals which represents the spatial dependence |

*Transmission*

| | |
|---|---|
| $S$ | $N \times K$ network incidence matrix, whose $k^{th}$ column contains $\pm 1$ for the rows corresponding to the terminal nodes (buses) of circuit k; and zero for the others |
| $f_{t,\tau}$ | $K$-dimensional vector of circuit flows $\{f_{t,\tau,k}\}$ |
| $e_{t,\tau}$ | $N$-dimensional vector of hydro generation. The energy production $e_{t,\tau,i}$ of each hydro $i$ is in the row of its respective network bus, $n(i)$ (all other values are zero). |
| $g_{t,\tau}$ | $N$-dimensional vector of thermal generation. The energy production $g_{t,\tau,j}$ of each thermal plant $j$ is in the row of its respective network bus, $n(j)$ (all other values are zero). |
| $\hat{d}_{t,\tau}$ | $N$-dimensional vector of load, where each power injection is in the row of its respective network bus (all other values $= 0$) |
| $r_{t,\tau}^{*s}$ | $N$-dimensional vector of renewable generation. For the existing renewable plants $(r \in \mathbb{R}_e), r_{t,\tau}^{*s} = \hat{r}_{t,\tau,r}^{s}$. For the candidate renewable plants $(r \in \mathbb{R}_x), r_{t,\tau}^{*s} = \hat{r}_{t,\tau,r}^{s} \times x_{t,r}^{*}$. |
| $\gamma_k$ | susceptance of circuit $k$ |
| $\theta_{t,\tau,F(k)}$ | nodal voltage angles at the "from" terminal bus of circuit $k$, represented as $F(k)$ |
| $\theta_{t,\tau,T(k)}$ | nodal voltage angles at the "to" terminal bus of circuit $k$, represented as $T(k)$ |
| $M_k$ | "big M" parameter |

## II. INTRODUCTION

The biggest challenge for the transmission planning in hydrothermal systems is the need to design a network that accommodates different hydro dispatch patterns (which, in turn, have a degree of flexibility and can be rearranged to accommodate the transmission constraints themselves). On the other hand, systems with increasing intermittent renewable energy penetration, more than presenting non-dispatchable generation, have an even higher dispatch variability, which in turn leads to the necessity of network robustness in order to meet the different dispatch scenarios. In consequence, representing uncertainties is a key issue and planning system's expansion in an economic efficient way is a not a trivial task.

Economies of scale is also an issue because generators are free to decide when and where to build new capacity, and in consequence planners have to take into account an additional degree of uncertainty when designing the network. In order to exemplify that, if three hydro plants totaling 4,500 MW are built in a given river basin located 2,000 km from the main grid, it may be more economic to use a transmission system with a higher voltage level than usual, e.g., 800 kV. On the other hand, if, in practice, only one of the plants is built, or the second one is installed years apart and the third not, the economy of scale is lost, and a sequence of 500 kV transmission systems could have been a better choice. The same issue may occur with wind parks due to the gradual penetration (because of the CAPEX reduction curves, market, auctions, etc.).

Another challenge is on the design of open-access transmission tariffs. The economic evaluation of new generation projects has to take into account an estimate of the associated transmission charges. The reason is that new hydro plants, which are usually farther from load centers - and thus have higher transmission costs - compete against new gas-fired thermal plants, which are closer to load centers, with correspondingly smaller associated transmission costs (if they in turn reflect the locational factors). More recently, both hydro and gas plants compete against cogeneration plants, which are also closely connected to the load centers and usually have even smaller transmission costs associated to transmission credits based on incentives and/or savings in transmission losses). Because transmission charges may vary depending on plant location, they interfere in the technology/plant competitiveness.

In other words, transmission planners have a "chicken or the egg" dilemma. On the one hand, generation investors have to know beforehand their associated transmission costs, in order to factor them in their contract prices. On the other hand, transmission planners have to know which generators are going to be build (i.e. have won supply contracts) in order to design the transmission system reinforcements that meet the different dispatch scenarios, and thus be able to allocate the transmission costs.

In summary, the planning process should be able to jointly consider all available alternatives and perform *trade-off* analysis between investment cost and operating cost of each alternative in search of more cost-effective solution. Additionally, it should jointly consider the expansion of the generation system and transmission network [1], [2], [3].

However, solving the generation and transmission expansion problem simultaneously might be computationally hard to solve because of the combinatorial nature of this problem, the size of real systems and the need for robustness in the face of stochastic renewables. To maintain a manageable model size, some tools proposed in the literature apply a horizon decomposition heuristic, for example considering annual investment stages, that is, a problem of co-optimization of the investment and operation is solved for each year in a rolling horizon scheme [3]. Since the number of constraints increase linearly with the number of the dispatch scenarios considered in the operation, other models propose scenario reduction framework to select representative scenarios to be incorporated in the investment module. As an example, in [4], after applying this technique, a Progressive Hedging (PH) algorithm is used to solve the resulting reduced stochastic optimization model. On the other hand, the convergence of PH to an optimal solution is not guaranteed for the mixed-integer linear optimization problem (MILP).

Other papers apply a hierarchical expansion planning procedure that consists of two solution steps, in the first step the Generation Expansion Plan (GEP) is found and in the second one, the transmission expansion plan is encountered taking the GEP and dispatch scenarios (obtained using production cost simulation tools) into account [5], [6], [7]. As can be seen, transmission reinforcements are obtained in order to accommodate the generation investment decisions that were made without representing the transmission network in detail.

Another simplification commonly found in methods that analyze the expansion of transmission systems is the relaxation of the temporal nature of the problem, i.e., several methods proposed in the literature for transmission network expansion are generally static and do not take into account the economies of scale between the stages among the study horizon [3], [4], [8], [9].

This paper aims at solving the "chicken or the egg" dilemma by finding an optimal generation and transmission (G&T) expansion plan. The G&T problem is formulated as an optimization problem and can be solved by a decomposition scheme based on a two-stage approach, as described below:

- **First-stage problem (the investment sub-problem):** formulated as a MILP problem where the objective is to propose alternatives for the G&T expansion plan;
- **Second-stage problem (the operation sub-problem):** the objective of the second stage is to evaluate the performance of the expansion alternatives proposed in the first-stage, producing the results that will be used in the first stage to improve the expansion solution. The second stage is solved by a probabilistic dispatch simulation tool.

In summary, the G&T expansion planning task is performed through a computational tool which determines the least-cost expansion plan for an electricity system dealing with hydro, thermal, variable renewable energy sources (VREs) and transmission candidate projects. The least-cost G&T plan is achieved by optimizing the *trade-off* between investment costs to build new projects and the expected value of operative costs obtained from the stochastic hydrothermal dispatch model, which allows a detailed representation of system's operation under uncertainty respecting network flows and limits. Finally, it is worth mentioning that this Benders decomposition scheme guarantees the optimal solution for this problem.

### III. GENERATION & TRANSMISSION EXPANSION PLAN

#### A. Overview of the Methodology

The diagram below illustrates the main features of the methodology.

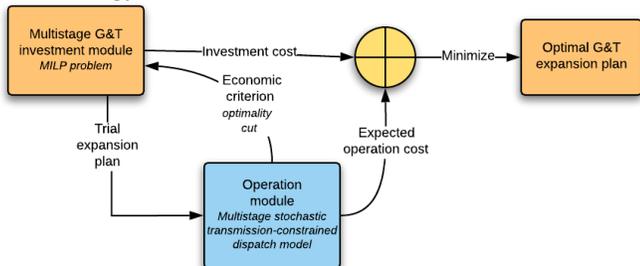

Fig. 1. Decomposition scheme for stochastic G&T expansion planning.

The Benders decomposition scheme was originally applied in the MODPIN planning model, has extensively been used in Latin America [10], [11], [12], [13]. In 1994, MODPIN was replaced by the OPTGEN model [14], [15].

The investment module is a MILP problem that produces a trial generation-transmission expansion plan. This trial plan is sent to the operation model, which carries out a multistage transmission-constrained stochastic optimization of system operation using the Stochastic Dual Dynamic Programming (SDDP) algorithm. Benders cuts representing the derivatives of expected operation cost with respect to the investment decisions are then calculated from SDDP's simulation of system operation for a large number of probabilistic scenarios (hydro inflows, VRE production, load uncertainty, equipment outages etc.).

#### B. Investment Module

We show below the investment problem after $M$ iterations.

$$Min \sum_t \sum_{i \in \mathbb{I}_x} I_i x_{t,i} + \sum_{j \in \mathbb{J}_x} I_j x_{t,j} + \sum_{r \in \mathbb{R}_x} I_r x_{t,r} + \sum_{k \in \mathbb{K}_x} I_k x_{t,k} + w \quad (1a)$$

$$x_{t,i} \geq x_{t-1,i} \qquad \forall t \in T; i \in \mathbb{I}_x \quad (1b)$$

$$x_{t,j} \geq x_{t-1,j} \qquad \forall t \in T; j \in \mathbb{J}_x \quad (1c)$$

$$x_{t,r} \geq x_{t-1,r} \qquad \forall t \in T; r \in \mathbb{R}_x \quad (1d)$$

$$x_{t,k} \geq x_{t-1,k} \qquad \forall t \in T; k \in \mathbb{K}_x \quad (1e)$$

$$w \geq \sum_t \left( \sum_{i \in \mathbb{I}_x} \bar{\mu}_{t,i}^m x_{t,i} + \sum_{j \in \mathbb{J}_x} \bar{\mu}_{t,j}^m x_{t,j} + \sum_{r \in \mathbb{R}_x} \bar{\mu}_{t,r}^m x_{t,r} + \sum_{k \in \mathbb{K}_x} \bar{\mu}_{t,k}^m x_{t,k} + \bar{\mu}_{0t}^m \right) \qquad \forall m \quad (1f)$$

*Mutually exclusive projects*

$$\sum_{p \in P_c^{ex}} \sum_t x_{t,p} \leq 1 \qquad \forall c \in R^{ex} \quad (1g)$$

*Associated projects*

$$\sum_t x_{t,p1} - \sum_t x_{t,p2} = 0 \qquad \forall p1, p2 \in P_c^{as}, \forall c \in R^{as} \quad (1h)$$

*Precedence between projects*

$$\sum_t x_{t,p2} - \sum_t x_{t,p1} \geq 0 \qquad \forall i, j \in P_c^{pre}, \forall c \in R^{pre} \quad (1i)$$

Note: for simplicity, we do not represent in this formulation additional real-life constraints, such as firm energy/capacity constraints, VRE penetration targets and others. For more info, the reader should refer to [15][15].

### C. Convergence criterion

The optimal solution value of the investment module in each iteration is a Lower Bound (LB) for the overall optimal solution because the linear representation of the operating cost is an under-approximation of the true cost. In turn, the sum of the investment cost of the trial expansion plan and the "real" operating cost (calculated by the operation module) is an Upper Bound (UB) for the global optimum, because this is a feasible solution, not necessarily the optimal one.

Because the incorporation of the Benders cuts in each iteration successively improves the operative cost function approximation in the investment module, the LB progressively increases. In turn, the UB progressively decreases, since it is only updated when a better solution is found. Therefore, we know that the global optimum has been achieved when Upper and Lower Bounds coincide (within a user-specified tolerance).

## IV. OPERATION MODULE

Given the trial optimal investment decisions $\{x_{t,i}^*\}$, $\{x_{t,j}^*\}$, $\{x_{t,r}^*\}$ and $\{x_{t,k}^*\}$ of the investment module in the $M$-th iteration of the Benders decomposition scheme, we solve the stochastic scheduling problem using the SDDP algorithm, described next.

The optimal operation problem for cost-based power systems is a classical and well-studied problem since the famous solution SDDP stated in the seminal paper by Pereira et al. [16] (which is an enhanced publication of the work presented in [17]), especially for the case of hydrothermal systems.

Because of its great practical success, the SDDP algorithm has been studied continuously. The method has deeply been analyzed in light of stochastic programming framework [18], the performance has been studied in [19], the method convergence has been analyzed by [18], [20] and the stopping criteria have been also studied in [21].

It is plausible to mention that many other methods and variants were proposed to solve the operation problem such as the approximate dynamic programming [22]. On the other hand, we solve the stochastic scheduling problem using the SDDP algorithm, because from the authors' perspective, the state of art solution of the long-term hydrothermal dispatch problem is still given by the SDDP algorithm and also because it is still on the leading frameworks to solve the problem and it is heavily used in industrial applications.

### A. SDDP Formulation

The figure below shows the main components of the operation problem for stage $t$, scenario $s$:
1. SDDP state variables at the beginning of the stage (in this example, initial storage $v(t)$ and hydro inflow along the stage, $a(t)$);
2. Reservoir storage balance equations, which determine the hydro turbined outflow, $u(t)$;
3. Power balance equation, which determines the least-cost operation of the thermal plants required to meet the residual load (after subtracting hydro generation and VRE production). In the SDDP formulation, the resulting operation cost is known as the immediate cost function (ICF);
4. Future cost functions (FCF) $l = 1, \dots, L$ of the SDDP state variables for the next stage: the final storage $v(t+1)$ and $l = 1, \dots, L$ conditioned inflow scenarios $a(t+1, l)$.
5. The objective function: minimize the sum of immediate cost (ICF) and the mean future cost ($1/L \sum FCF_l$).

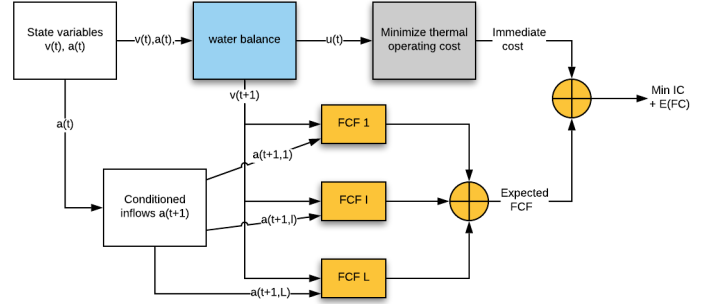

Fig. 2. Main components of SDDP's operation problem for stage $t$, scenario $s$.

*Objective function*

$$\alpha_t(\hat{v}_t^s, \hat{a}_t^s) = Min \sum_j c_j \sum_\tau g_{t,\tau,j} + \frac{1}{L}\sum_l \alpha_{t+1}^l \quad (2a)$$

*Storage balance for each stage*

$$v_{t+1,i} = \hat{v}_{t,i}^s + \hat{a}_{t,i}^s - (u_{t,i}+v_{t,i}) + \sum_{u \in \mathcal{U}_i}(u_{t,u}+v_{t,u}) \quad \forall i \quad \leftarrow \pi_{t,i}^h \quad (2b)$$

Note: for notational simplicity, we do not represent in this formulation real-life features of the storage balance equations such as evaporation, filtration, water diversion for irrigation and city supply, transposition and others.

*Storage limits*

For the existing hydro plants, the storage limit is a given value, $\overline{v}_i$

$$v_{t+1,i} \leq \overline{v}_i \qquad \forall i \in \mathbb{I}_e \quad (2c)$$

For the candidate hydro plants, the storage limit $\overline{v}_i^*$ depends on the investment decision $x_{t,i}^*$ at the current iteration of the Benders decomposition scheme (investment module):

$$v_{t+1,i} \leq \overline{v}_i^* (= \overline{v}_i \times x_{t,i}^*) \qquad \forall i \in \mathbb{I}_x \quad \leftarrow \pi_{t,i}^{\overline{v}} \quad (2d)$$

Note that the SDDP operating module "sees" both existing and candidate storage limits $\overline{v}_i$ and $\overline{v}_i^*$ as given values. In other words, the operating module does not "know" that is being run as part of a Benders decomposition scheme with the investment module. The investment information is only used explicitly in the calculation of the Benders cuts from the operation to the investment module. One advantage of this scheme that the same model used in operations studies can be used in planning studies, without any modification.

*Turbined outflow limits*

$$u_{t,i} \leq \overline{u}_i \qquad \forall i \in \mathbb{I}_e \qquad (2e)$$

$$u_{t,i} \leq \overline{u}_i^*(= \overline{u}_i \times x_{t,i}^*) \qquad \forall i \in \mathbb{I}_x \quad \leftarrow \pi_{t,i}^{\overline{u}} \qquad (2f)$$

The same logic of the storage limit $(2d)$ applies to the turbined outflow $(2f)$.

*Hydro generation*

$$e_{t,i} = \phi_i u_{t,i} \qquad \forall i \qquad (2g)$$

For notational simplicity, hydro generation is represented here as a linear function of the turbined outflow. In real-life, other factors are modeled such as the variation of reservoir head with storage, increase of tailwater with total outflow and encroachment effect associated to downstream plants.

$$\sum_\tau e_{t,\tau,i} = e_{t,i} \qquad \forall i \qquad (2h)$$

$$e_{t,\tau,i} \leq \overline{e}_i \qquad \forall i \qquad (2i)$$

Note that it is not necessary to have investment variables for hydro energy production because this is already done for turbined outflow.

*Thermal generation*

As in the hydro case, the generation capacity of candidate plants changes with the investment module iterations.

$$g_{t,\tau,j} \leq \overline{g}_j \qquad \forall j \in \mathbb{J}_e \qquad (2j)$$

$$g_{t,\tau,j} \leq \overline{g}_j^*\left(= \overline{g}_j \times x_{t,j}^*\right) \qquad \forall j \in \mathbb{J}_x \quad \leftarrow \pi_{t,\tau,j}^{\overline{g}} \qquad (2k)$$

As in the previous cases, we show a simple representation of thermal generation in this formulation. In actual applications, there are models for efficiency curves, multiple fuels, fuel storage, fuel contracts and unit commitment.

*Renewable generation*

Renewable generation is represented as energy production scenarios $\{\hat{r}_{t,\tau,r}^s\}$ in the power balance equations, described next.

*Transmission network equations and constraints*

The first set of equations represents the power balance in each bus (Kirchhoff's first law):

$$Sf_{t,\tau} + e_{t,\tau} + g_{t,\tau} = \hat{d}_{t,\tau} - r_{t,\tau}^{*s} \qquad \leftarrow \pi_{t,\tau,j}^d \qquad (2l)$$

Next, we represent Kirchhoff's second law:

$$f_{t,\tau,k} = \gamma_k\big(\theta_{t,\tau,F(k)} - \theta_{t,\tau,T(k)}\big) \qquad \forall k \in \mathbb{K}_e \qquad (2m)$$

For the candidate circuits, the second law is represented as the following constraint:

$$\big|f_{t,\tau,k} - \gamma_k(\theta_{t,\tau,F(k)} - \theta_{t,\tau,T(k)})\big| \leq \Delta_{t,k}^*\big(= M_k[1 - x_{t,k}^*]\big) \qquad \leftarrow \pi_{t,\tau,k}^\gamma \qquad (2n)$$

We can see that if the candidate circuit is built in the current Benders iteration ($x_{t,k}^* = 1$), constraint $(2n)$ becomes equal to equation $(2m)$ of the existing circuits. Conversely, if $x_{t,k}^* = 0$, constraint $(2n)$ is relaxed.

It is plausible to mention that if $M_k$ is arbitrarily big, the mathematical optimization problem becomes ill-conditioned. Therefore, for each candidate right-of-way the smallest value of $M_k$ capable of enforcing (or relaxing) the constraint when needed is calculated. Initially suppose that there is an existent circuit having reactance $\gamma_k^0$, capacity $f_k^0$ and the same bus terminals as candidate circuit $k$. The maximum angle difference between these bus terminals is $f_k^0/\gamma_k^0$; therefore one can set $M_k = \gamma_k (f_k^0/\gamma_k^0)$. For a new corridor connecting buses $i_k$ and $j_k$, i.e., no existing circuit connect the bus terminals, the maximum angle difference can be derived considering each path from $i_k$ to $j_k$ composed by existing circuits. For each such circuit, its maximum angle difference is the ratio mentioned earlier, and summing these terms results in the maximum angle difference between $i_k$ and $j_k$. Since there may be several paths connecting buses $i_k$ and $j_k$, the smallest value of $M_k$ will be the candidate's reactance times the length of the shortest path between $i_k$ and $j_k$ (a circuit "length" is the ratio of its capacity and its reactance) [23], [24].

The use of the aforementioned disjunctive formulation to solve benchmark problems found in the transmission expansion literature was proved to be very effective, they were solved faster and the optimal solution was obtained and proven as detailed in [24].

Finally, the circuit flow limits are represented as:

$$|f_{t,\tau,k}| \leq \overline{f}_k \qquad \forall k \in \mathbb{K}_e \qquad (2o)$$

$$|f_{t,\tau,k}| \leq \overline{f}_{t,k}^*(= \overline{f}_k x_{t,k}^*) \qquad \forall k \in \mathbb{K}_x \quad \leftarrow \pi_{t,\tau,k}^{\overline{f}} \qquad (2p)$$

*Conditioned inflow scenarios for t+1*

For simplicity of presentation, we show a multivariate $AR(1)$ model. In practice, SDDP uses a multivariate periodic autoregressive $(PAR(p))$ model with up to six past time stages:

$$\frac{(a_{t+1,i}^l - \hat{\mu}_{t+1,i})}{\hat{\sigma}_{t+1,i}} = \hat{\rho}_{t,i} \times \frac{(\hat{a}_{t,i}^s - \hat{\mu}_{t,i})}{\hat{\sigma}_{t,i}} + \sqrt{1 - \hat{\rho}_{t,i}^2} \times \hat{\xi}_{t,i}^l \quad \forall i \leftarrow \pi_{t,i}^a \quad (2q)$$

Note: For clarity of presentation, the stochastic streamflow models are shown explicitly. In the actual SDDP implementation, they are represented implicitly.

*Future Cost Functions (FCFs)*

As it is well known in SDDP, the FCFs are represented by a set of hyperplanes:

$$\alpha_{t+1}^l \geq \sum_i \hat{\varphi}_{t+1,i}^{hp} \times v_{t+1,i} + \sum_i \hat{\varphi}_{t+1,i}^{ap} \times a_{t+1,i}^l + \hat{\varphi}_{t+1}^{op} \quad \forall p, l \quad (2r)$$

### B. Benders Cut to the Investment Module

As it is also well known, the SDDP algorithm is composed of three steps: (i) backward recursion; (ii) forward simulation; and (iii) convergence check. Here, we describe a fourth step used in planning models, which is the calculation of marginal capacity information for a new $(M+1)^{th}$ Benders cut to the investment module [12]:

$$w \geq \sum_t \left( \sum_{i \in \mathbb{I}_x} \bar{\mu}_{t,i}^{M+1} x_{t,i} + \sum_{j \in \mathbb{J}_x} \bar{\mu}_{t,j}^{M+1} x_{t,j} + \sum_{r \in \mathbb{R}_x} \bar{\mu}_{t,r}^{M+1} x_{t,r} + \sum_{k \in \mathbb{K}_x} \bar{\mu}_{t,k}^{M+1} x_{t,k} + \bar{\mu}_{0t}^{M+1} \right) \quad (3)$$

The Benders cut coefficients are obtained from the multipliers associated to the following constraints: (i) hydro storage and turbined outflow limits $(2d)$ and $(2f)$; (ii) thermal generation capacity limits $(2k)$; (iii) energy production for the VREs $(2l)$; and (iv) Kirchhoff's second law and flow limits for the transmission components $(2n)$ and $(2p)$. The calculation of these coefficients is presented below:

$$\bar{\mu}_{t,i}^{M+1} = \frac{1}{S} \sum_s \left( \bar{v}_i \times \pi_{t,i}^{\bar{v}s} + \bar{u}_i \times \pi_{t,i}^{\bar{u}s} \right) \quad (5a)$$

$$\bar{\mu}_{t,j}^{M+1} = \frac{1}{S} \sum_s \sum_\tau \left( \bar{g}_j \times \pi_{t,\tau,j}^{\bar{g}s} \right) \quad (5b)$$

$$\bar{\mu}_{t,r}^{M+1} = \frac{1}{S} \sum_s \sum_\tau \left( \hat{r}_{t,\tau,r}^s \times \pi_{t,\tau,n(r)}^{ds} \right) \quad (5c)$$

$$\bar{\mu}_{t,k}^{M+1} = \frac{1}{S} \sum_s \sum_\tau \left( -M_k \times \pi_{t,\tau,k}^\gamma + \bar{f}_k \times \pi_{t,\tau,k}^{\bar{f}} \right) \quad (5d)$$

$$\bar{\mu}_{0t}^{M+1} = \frac{1}{S} \sum_s \sum_\tau \sum_j c_j g_{t,\tau,j}^{*s} - \sum_{i \in \mathbb{I}_x} \bar{\mu}_{t,i}^{M+1} x_{t,i}^* - \sum_{j \in \mathbb{J}_x} \bar{\mu}_{gt,j}^{M+1} x_{t,j}^* - \sum_{r \in \mathbb{R}_x} \bar{\mu}_{t,r}^{M+1} x_{t,r}^* - \sum_{k \in \mathbb{K}_x} \bar{\mu}_{t,k}^{M+1} x_{t,k}^* \quad (5e)$$

## V. CASE STUDY

### A. The Bolivian System

The integrated G&T planning approach proposed in this paper is applied to a case study based on the Bolivia power system. We assume that the optimization model can make no expansion decisions during the 2019-2021 sub-horizon and therefore, our case study starts at January 2022. In the beginning of the study horizon, there are 48 hydroelectric power plants (HPPs) presenting 923 MW of installed capacity, 98 thermal power plants (TPPs) with approximately 2,400 MW of installed capacity and 276 MW of VREs. The high voltage network is composed by 169 transmission lines (3 being of 24 kV, 34 of 69 kV, 85 of 115 kV and 47 of 230 kV) and 54 transformers. The figure presented below illustrates the aforementioned system:

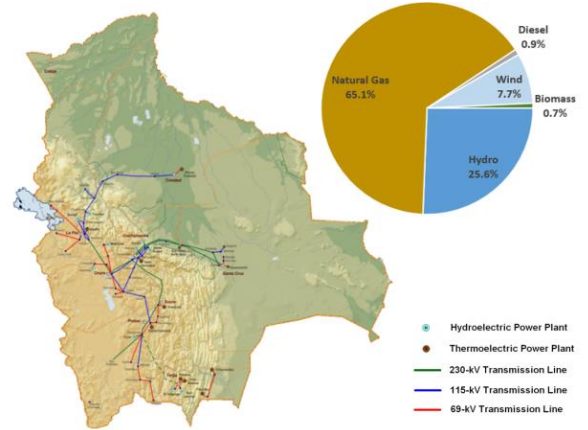

Fig. 3. The Bolivian system general overview.

The time intervals used in this study are monthly stages with load duration curve modeling and precise numbers of hours in each of the five blocks in each month (according to the month duration in order to maintain the following percentage of block durations from one to five, respectively: 2.3%, 3.6%, 9.6%, 48.2% and 36.3%. In 2022, a demand of 13,805 GWh demand is expected with an average growth of 4.3% per year within the study horizon.

Regarding the uncertainty representation, 32 forward and backward scenarios were used in this case study, where each one is a combination of inflows for the hydro plants and generation scenarios for VREs maintaining the temporal and spatial correlations). Furthermore, it is plausible to mention that the study horizon taken into account in this paper is 2022-2030 and decisions on building (or not) new projects start at 2022.

## B. Candidate Projects and G&T Expansion Assessment

In order to expand the generation matrix, the following candidate projects were taken into account: (i) 9 combined cycle gas turbine plants (CCGTs) of 250 MW; (ii) 9 open cycle gas turbine plants (OCGTs) of 100 MW; (iii) 9 diesel plants of 80 MW; (iv) 6 GW of wind; and (v) 10 GW of solar projects. The investment data considered in this case study is summarized in the table presented below:

TABLE I
INVESTMENT DATA OF THE GENERATION PROJECTS

|  | CCGT | OCGT | Diesel | Solar | Wind |
|---|---|---|---|---|---|
| Investment cost ($/kW) | 900 | 700 | 700 | 1200 | 1400 |
| Payments during construction (%) | 33-33-33 | 50-50 | 50-50 | 50-50 | 50-50 |
| Lifetime (years) | 20 | 20 | 20 | 25 | 25 |
| Fixed O&M cost ($/kW.year) | 25 | 15.3 | 12 | 20 | 20 |
| WACC (% p.a.) | 9% | 9% | 9% | 9% | 9% |

Additionally, for the solar and wind projects we have considered a CAPEX reduction curve of 3% p.a. and 1% p.a., respectively. Regarding network alternatives, 94 transmission circuits were considered as candidate projects (7 being transmission lines of 69 kV, 31 of 115 kV, 30 of 230 kV and finally 26 transformers). It is plausible to emphasize that all candidate projects present real electrical parameters, costs and lifetime with respect to the Bolivian system and are represented as binary decision variables in the investment module.

The next step is to determine the optimal G&T expansion plan. As explained in previous sections, the applied methodology is basically a Benders decomposition scheme. In summary, the investment module will solve the new capacity problem by evaluating the trade-off between investing in each project and the associated impact on the operating costs. At each iteration, a new investment/operation solution is found and in consequence, the total cost of the best solution found so far (the so called Upper Bound) and the total *expected* cost (the so called Lower bound) are updated. Every time a cheaper feasible solution is found, the UB is updated.

Furthermore, besides determining *which* projects should be built, the investment module also decides *when* the projects should come into operation along the study horizon. As can be seen, the CAPEX-OPEX trade-off is optimized by solving the *sizing* and *timing* problems simultaneously in order to find the least-cost expansion plan.

The expansion planning problem of this case study was run on a virtual server on AWS with 36 CPUs, each CPU being a hardware hyperthread on a 3.0 GHz Intel Xeon Platinum 8000-series processor and 72 GB of RAM. The numbers of Benders iterations (solving the investment problem MILP plus a SDDP-based stochastic operation model) was 80 to reach the target gap of 3% in 151 minutes. The convergence process is illustrated in the figure presented below:

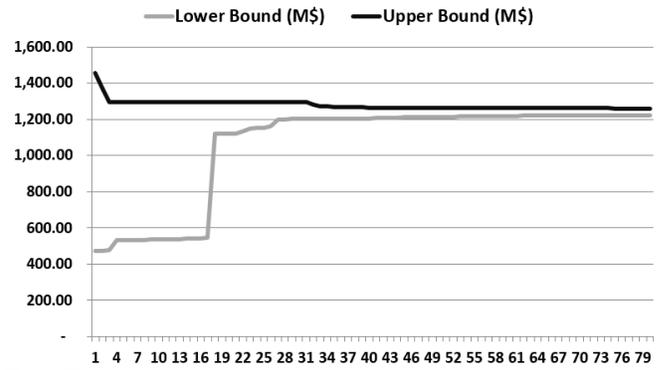
Fig. 4. The convergence process.

Besides the convergence process, Figure 5 shows the added generation capacity per technology in each year; Figure 6 presents the consequent installed capacity versus peak load comparison among the study horizon; and finally, Table II shows the transmission circuits which are part of the optimal G&T expansion plan.

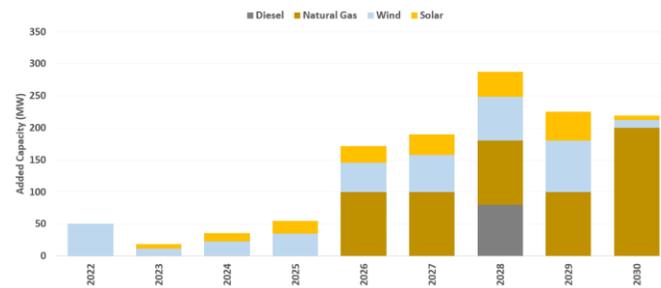
Fig. 5. Added generation capacity per technology in each year.

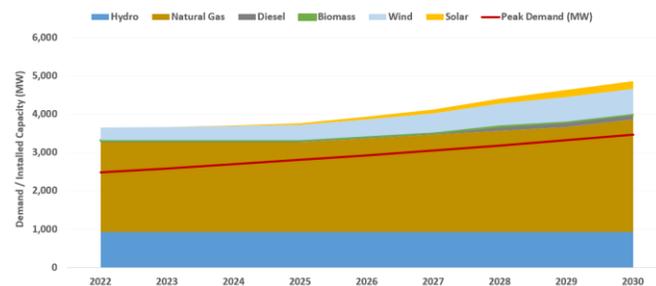
Fig. 6. Installed capacity *vs* peak demand after the optimal G&T expansion plan is found

TABLE II
TRANSMISSION CIRCUITS PART OF THE OPTIMAL G&T EXPANSION PLAN

| Bus From Name | Bus From Voltage (kV) | Bus To Name | Bus To Voltage (kV) | Circuit Cost (M$) | Lifetime (Years) | Circuit Rating (MW) | Circuit Type | Entry Year |
|---|---|---|---|---|---|---|---|---|
| SUC-230 | 230 | SUC-069 | 69 | 2.89 | 30 | 57 | Transformer | 2022 |
| VIN-230 | 230 | VIN-115 | 115 | 4.71 | 30 | 95 | Transformer | 2025 |
| ARJ-069 | 69 | SUC-069 | 69 | 1.62 | 30 | 42 | Line | 2028 |
| WAR-230 | 230 | BEL-230 | 230 | 6.46 | 30 | 275 | Line | 2029 |
| VIN-115 | 115 | VIN-069 | 69 | 2.21 | 30 | 47.5 | Transformer | 2029 |
| WAR-230 | 230 | URU-230 | 230 | 11.96 | 30 | 275 | Line | 2029 |
| BRE-230 | 230 | BRE-069 | 69 | 4.11 | 30 | 142.5 | Transformer | 2029 |
| GCH-069 | 69 | PIN-069 | 69 | 1.27 | 30 | 95 | Line | 2029 |
| URU-230 | 230 | URU-069 | 69 | 4.11 | 30 | 142.5 | Transformer | 2029 |
| URU-230 | 230 | URU-069 | 69 | 4.11 | 30 | 142.5 | Transformer | 2029 |
| GCH-069 | 69 | PAR-069 | 69 | 1.04 | 30 | 128 | Line | 2029 |
| GCH-069 | 69 | ZOO-069 | 69 | 1.24 | 30 | 89 | Line | 2029 |
| VHE-115 | 115 | IRP-115 | 115 | 8.36 | 30 | 74 | Line | 2030 |
| LIT-230 | 230 | LIT-115 | 115 | 4.56 | 30 | 71 | Transformer | 2030 |
| CAT-115 | 115 | CAT-069 | 69 | 1.95 | 30 | 23.75 | Transformer | 2030 |
| CAT-115 | 115 | CAT-069 | 69 | 1.95 | 30 | 23.75 | Transformer | 2030 |

It is worth mentioning that due to the VRE penetration, only OCGTs have been selected by the model with respect to the natural gas alternatives. Since OCGTs have lower CAPEX, however, higher OPEX, in general they are built to meet the peak demand (plants also known as *peakers*) and operational situations where VREs are not generating, i.e., they are built by the fact of not dispatching frequently, otherwise the break-even would be achieved and the CCGTs would be more cost-effective as they are more efficient.

The optimal G&T expansion plan avoids energy deficits, overloading of transmission circuits and load shedding at all buses. As the CAPEX-OPEX trade-off analysis has been successfully performed, the final SDDP simulation taking the optimal G&T expansion plan into account results in stable short-run marginal costs, around 15 $/MWh in the long-term period.

## VI. Conclusions

In this paper, the authors describe an integrated generation-transmission expansion planning methodology based on Benders decomposition and multistage stochastic optimization. In that way, the G&T expansion problem is formulated as an optimization problem and can be solved by a decomposition scheme. It can be seen that a real system planning process may be represented in detail (i) in the investment module and also (ii) in the operating module. Taking the high level of the problem's complexity mainly due to the size of the system and the combinatorial nature of the expansion planning problem (sizing and timing), this paper has also proven that using the virtual server on AWS with 36 CPUs, the planning problem of a real system can be solved in a low computational time achieving the optimal G&T expansion plan.